\documentclass[english]{amsart}
\usepackage[T1]{fontenc}
\usepackage[latin9]{inputenc}
\usepackage{amsthm}
\usepackage{graphicx}
\usepackage{amssymb}

\numberwithin{equation}{section} %% Comment out for sequentially-numbered
\numberwithin{figure}{section} %% Comment out for sequentially-numbered
\theoremstyle{plain}
\theoremstyle{plain}
\newtheorem{thm}{Theorem}
  \theoremstyle{remark}
  \newtheorem{rem}[thm]{Remark}
  \theoremstyle{plain}
  \newtheorem{lem}[thm]{Lemma}

\usepackage{babel}

\begin{document}

\title[Interpolation of compact Lipschitz operators - Counterexamples]{Counterexamples for interpolation of compact lipschitz operators}

\author{Michael Cwikel and Alon Ivtsan}

\address{Department of Mathematics, Technion - Israel Institute of Technology,
Haifa 32000, Israel}

\email{mcwikel@math.technion.ac.il \textit{and} aloniv@techunix.technion.ac.il}

\thanks{The research of the first named author was supported by the Technion
V.P.R.\ Fund and by the Fund for Promotion of Research at the Technion. }


\begin{abstract}
Let $\left(A_{0},A_{1}\right)$ and $\left(B_{0},B_{1}\right)$ be
Banach couples with $A_0\subset A_1$ and $B_0 \subset B_1$ and let $T:A_1 \to B_1$ be a
possibly nonlinear compact Lipschitz map whose restriction to $A_0$ is
also a compact Lipschitz map into $B_0$.
It is known that $T$ maps $\left(A_{0},A_{1}\right)_{\theta,q}$ boundedly
into $\left(B_{0},B_{1}\right)_{\theta,q}$
for each $\theta\in(0,1)$ and $q\in[1,\infty]$ and that this map
is also compact if $T$ is linear. We present examples which show
that in general the map $T:\left(A_{0},A_{1}\right)_{\theta,q}
\to\left(B_{0},B_{1}\right)_{\theta,q}$
is not compact.
\end{abstract}
\maketitle

\section{Introduction}

Let us begin by stating a theorem which was obtained in the the 1990's
in \cite{CKS} and \cite{Cw2}.
\begin{thm}
\label{thm:ccks}Let $\left(A_{0},A_{1}\right)$ and $\left(B_{0},B_{1}\right)$
be Banach couples. Suppose that $T:A_{0}+A_{1}\to B_{0}+B_{1}$ is
a linear operator which maps $A_{0}$ to $B_{0}$ compactly, and $A_{1}$
to $B_{1}$ boundedly. Then $T$ maps the Lions--Peetre space $\left(A_{0},A_{1}\right)_{\theta,q}$
to $\left(B_{0},B_{1}\right)_{\theta,q}$ compactly for each $\theta\in(0,1)$
and each $q\in[1,\infty]$.
\end{thm}
Various special cases of Theorem \ref{thm:ccks} go back to the 1960's.
Mark Krasnosel'skii \cite{kras} gave the initial proof in the case
where all of $A_{0}$, $A_{1}$, $B_{0}$, $B_{1}$, $\left(A_{0},A_{1}\right)_{\theta,q}$
and $\left(B_{0},B_{1}\right)_{\theta,q}$ are $L^{p}$ spaces (of
course for possibly different values of the exponents $p$). Jacques-Louis
Lions and Jaak Peetre (see Théorème (2.1) and Théorème (2.2) of \cite{LP}
pp.~36--38) proved it in the case where $A_{0}=A_{1}$ and in the
case where $B_{0}=B_{1}$. Arne Persson \cite{P} proved it in the
case where the couple $\left(B_{0},B_{1}\right)$ satisfies a certain
{}``approximation hypothesis'' (see \cite{P} p.~216). K.~Hayakawa
\cite{hayakawa} proved it in the case where $T$ satisfies the additional
condition that $T:A_{1}\to B_{1}$ is also compact.

In this note we investigate the question of whether Theorem \ref{thm:ccks}
can be extended to cases where the operator $T$ is nonlinear. This
question seems natural since it has been possible to extend a considerable
part of the theory of Lions--Peetre interpolation spaces to the context
of nonlinear operators, in particular those operators which satisfy
appropriate Lipschitz conditions and boundedness conditions. This
has been done by Jaak Peetre in \cite{PeetreCluj} and by Jacques-Louis
Lions \cite{lions} and in rather more detail by Luc Tartar \cite{tartar}.
The papers \cite{lions} and \cite{tartar} also include some applications
of their nonlinear interpolation results to partial differential equations.
We are grateful to Lavi Karp for drawing our attention to the more
recent book \cite{RunstSickel} of Thomas Runst and Winfried Sickel
which includes a summary of results of this kind on pp.~87--92. We
also refer to \cite{browder} for another approach to extending results
about interpolation of linear operators to interpolation of Lipschitz
operators.

Here is a particularly simple instance of the kinds of results about
nonlinear operators which are presented in \cite{lions,PeetreCluj,RunstSickel,tartar}.
\begin{thm}
\label{thm:simple}Let $\left(A_{0},A_{1}\right)$ and $\left(B_{0},B_{1}\right)$
be Banach couples. Suppose that $A_{0}\subset A_{1}$ and
$B_{0}\subset B_{1}$. Let $T$ be
a (possibly nonlinear) map of $A_{1}$ into $B_{1}$ which satisfies
the following two properties:\begin{equation}
T\left(A_{0}\right)\subset B_{0}\mbox{ and }\left\Vert T(a)\right\Vert _{B_{0}}\le C_{0}\left\Vert a\right\Vert _{A_{0}}\,\mbox{for each }a\in A_{0}\,,\label{eq:sone}\end{equation}
and\begin{equation}
\left\Vert T(a)-T(a')\right\Vert _{B_{1}}\le C_{1}\left\Vert a-a'\right\Vert _{A_{1}}\,\mbox{for all }a,a'\in A_{1}\,.\label{eq:stwo}\end{equation}
where $C_{0}$ and $C_{1}$ are positive constants.

Then $T$ maps the space $\left(A_{0},A_{1}\right)_{\theta,q}$ boundedly
into $\left(B_{0},B_{1}\right)_{\theta,q}$ for each $\theta\in(0,1)$
and $q\in[1,\infty]$, and satisfies the estimate \[
\left\Vert T(a)\right\Vert _{\left(B_{0},B_{1}\right)_{\theta,p}}\le C_{0}^{1-\theta}C_{1}^{\theta}\left\Vert a\right\Vert _{\left(A_{0},A_{1}\right)_{\theta,p}}\,\mbox{for all }a\in\left(A_{0},A_{1}\right)_{\theta,p}\,.\]

\end{thm}
The proof of this theorem is an immediate consequence of some simple
calculations with $K$-functionals. See \cite{lions,PeetreCluj,tartar}.
For a similar result, where the condition $A_{0}\subset A_{1}$ is
not imposed, but instead $T$ is required to be a Lipschitz map also
from $A_{1}$ into $B_{1}$ and $q$ is finite, see Theorem 4.1 on
p.~278 of \cite{Cobos-LPLemma}.

In this note we take Theorem \ref{thm:simple} as our point of departure
and ask the following question: Suppose that $\left(A_{0},A_{1}\right)$,
$\left(B_{0},B_{1}\right)$ and $T$ satisfy all the hypotheses of
the theorem, and one extra condition, namely that $T$ maps $A_{0}$
into $B_{0}$ compactly, or, alternatively, that $T$ maps $A_{1}$
into $B_{1}$ compactly. Is either one of these extra conditions sufficient
to ensure that the bounded map by $T$ of $\left(A_{0},A_{1}\right)_{\theta,p}$
into $\left(B_{0},B_{1}\right)_{\theta,p}$ is also a compact map?

There are two special cases studied by Fernando Cobos \cite{Cobos-LPLemma},
which we will describe in a moment, where the answer to this question
is affirmative. However we shall see that, in general, the answer
to this question is negative. Furthermore the answer remains negative
even when we try imposing any or even all of the various above-mentioned
extra conditions which enabled Krasnosel'skii, Persson and Hayakawa
in turn to each prove their versions of Theorem \ref{thm:ccks} for
linear $T$. Nor does it help to also replace (\ref{eq:sone}) by
the apparently (see Remark \ref{rem:zzero}) stronger Lipschitz condition:
\begin{equation}
T\left(A_{0}\right)\subset B_{0}\,\mbox{and }\left\Vert T(a)-T(a')\right\Vert _{B_{0}}\le C_{0}\left\Vert a-a'\right\Vert _{A_{0}}\,\mbox{for all }a,a'\in A_{0}\,.\label{eq:sonebis}\end{equation}

In contrast to all these negative results, the above mentioned positive
results of Cobos show that it \textit{does} help to impose either
one of the extra conditions $A_{0}=A_{1}$ and $B_{0}=B_{1}$. I.e.,
he deals with the {}``nonlinear'' versions of each of the two cases
treated by Lions and Peetre \cite{LP}. In Cobos' results (see Theorem
2.1 on p. 274 of \cite{Cobos-LPLemma}) the condition (\ref{eq:sone})
has to be replaced by the Lipschitz condition (\ref{eq:sonebis}).
But he does not need to require that $A_{0}\subset A_{1}$ or
$B_{0}\subset B_{1}$. This means
that the condition $T:A_{1}\to B_{1}$ of Theorem \ref{thm:simple}
has to be reformulated and in fact replaced by the two conditions
$T:A_{0}+A_{1}\to B_{0}+B_{1}$ and $T\left(A_{1}\right)\subset B_{1}$.
Cobos shows that if the map $T:A_{j}\to B_{j}$ is compact for at
least one of the two values $j=0$ and $j=1$, then this suffices
to ensure the compactness of
$T:\left(A_{0},A_{1}\right)_{\theta,q}\to\left(B_{0},B_{1}\right)_{\theta,q}$
for each $\theta\in(0,1)$ and each $q\in[1,\infty)$. Of course,
in the cases that he is considering,
one has either $\left(A_{0},A_{1}\right)_{\theta,q}=A_{0}=A_{1}$
or $\left(B_{0},B_{1}\right)_{\theta,q}=B_{0}=B_{1}$.
\begin{rem}
\label{rem:zzero}In the case where $T$ maps the zero element $0$
of $A_{0}+A_{1}$ to the zero element of $B_{0}+B_{1}$ then (\ref{eq:sonebis})
is indeed a stronger condition than (\ref{eq:sone}). If $T(0)$ is
not the zero element, then, as in \cite{BrauerKarp}, we can consider
the auxiliary operator $\widetilde{T}$ defined by $\widetilde{T}(f)=T(f)-T(0)$.
Since $T(0)\in B_{0}\cap B_{1}$, the mapping properties and Lipschitz
properties of $T$ and of $\widetilde{T}$ are essentially equivalent
and (\ref{eq:sonebis}) for $T$ of course implies that $\widetilde{T}$
satisfies (\ref{eq:sone}).
\end{rem}
We shall present a counterexample, an example of a particular operator
$T$, which provides a negative answer to our question and also to
the other variants of that question mentioned above where one tries
to {}``save'' the situation by imposing extra conditions. In our
example the couples $\left(A_{0},A_{1}\right)$ and $\left(B_{0},B_{1}\right)$
will be one and the same. In fact we will have $A_{0}=B_{0}=L^{\infty}$
and $A_{1}=B_{1}=L^{1}$ where the underlying measure space is $[0,1]$
equipped with Lebesgue measure. Our operator $T$ will have the following
five properties.

\smallskip{}

\noindent \begin{flushleft}
$\begin{array}{ll}
\mbox{[i]} & T(0)=0,\ \mbox{where }0\mbox{ denotes the zero element of }A_{0}+A_{1}.\\
\\\mbox{[ii]} & T(A_{j})\subset A_{j}\ \mbox{and }\left\Vert T(a)-T(a')\right\Vert _{A_{j}}\le\left\Vert a-a'\right\Vert _{A_{j}}\,\mbox{for all }a,a'\in A_{j}\\
 & \mbox{and for }j=0,1\,.\\
\\\mbox{[iii}_{0}\mbox{]} & T\mbox{\ maps every bounded subset of }A_{0}\mbox{ into a relatively compact subset of }A_{0}.\\
\\\mbox{[iii}_{1}\mbox{]} & T\ \mbox{maps every bounded subset of }A_{1}\mbox{ into a relatively compact subset of }A_{1}.\\
\\\mbox{[iv]} & \mbox{For every }\theta\in(0,1)\ \mbox{and }q\in[1,\infty],\mbox{ the map }T:\left(A_{0},A_{1}\right)_{\theta,q}\to\left(A_{0},A_{1}\right)_{\theta,q}\\
 & \mbox{is not compact.}\end{array}$\smallskip{}

\par\end{flushleft}

We shall obtain this example in several steps. In Section \ref{sec:prelim}
we will collect some preliminary results. Then, in Section \ref{sec:onesided},
we will describe an operator, which we will denote by $T_{1}$, which
is a {}``one sided'' example, i.e., it has all of the above five
properties except {[}iii$_{0}${]}. Then the second {}``one sided''
example, an operator to be denoted by $T_{2}$, which will be presented
in Section \ref{sec:smallonesided}, will have all the above properties
except {[}iii$_{1}${]}. Finally in Section \ref{sec:twosided}, we
will see that the operator $T_{3}=T_{2}\circ T_{1}$, i.e., the composition
of our previous two examples, can serve as the promised {}``two sided''
example of an operator $T$ having all the above five properties.

In an earlier stage of this research we also obtained three other
examples, one of them considerably more elaborate than those of Sections
\ref{sec:onesided}, \ref{sec:smallonesided} and \ref{sec:twosided}.
Although it subsequently turned out that we can answer our particular
questions here without using these additional examples, we put them
on record in an appendix (Section \ref{sec:appendix}) in case they,
and/or the methods used for their construction, may ultimately prove
to be relevant for investigating other questions about interpolation
of Lipschitz operators.
\begin{rem}
\label{rem:comp}On several occasions we will use the obvious fact
that if any two operators both map $A_{0}+A_{1}$ into $A_{0}+A_{1}$
and satisfy conditions {[}i{]} and {[}ii{]}, then so does their composition.
\end{rem}

\begin{rem}
The operator $T$ in every one of our examples will have properties
{[}i{]} and {[}ii{]} and therefore also the boundedness property $\Vert T(a)\Vert_{A_{j}}\le\Vert a\Vert_{A_{j}}$
for each $a\in A_{j}$. The couple $\left(A_{0},A_{1}\right)$ in
all our examples always satisfies $A_{0}\subset A_{1}$, but this
is only a convenience rather than a necessity for their construction,
and it also makes the comparison with Theorem \ref{thm:simple} more
explicit. It is a trivial matter to obtain modified versions of our
examples where neither of the inclusions $A_{0}\subset A_{1}$ and
$A_{1}\subset A_{0}$ hold.
\end{rem}

\begin{rem}
It will be apparent that in every example presented in this paper,
including in the appendix, the range couple $\left(B_{0},B_{1}\right)$
satisfies Arne Persson's approximation hypothesis. This is because
in each case we have $\left(A_{0},A_{1}\right)=\left(B_{0},B_{1}\right)$
and $\left(A_{0},A_{1}\right)$ is either $(L^{\infty},L^{1})$ or
(in just one example) $\left(\ell^{1},\ell^{\infty}\right)$, and
so one can invoke the proposition on pp.~218--219 of \cite{P}. Furthermore,
in all examples of this paper, among the spaces $\left(A_{0},A_{1}\right)_{\theta,q}$
on which we will show that $T$ does not act compactly, will be either
$\ell^{p}$ or $L^{p}$ for some $p\in(1,\infty)$. Thus, in all of
our examples, we are using the same spaces as appear in Krasnosel'skii's
theorem.
\end{rem}

\begin{rem}
For convenience and flexibility in applications, it is natural to
formulate theorems like Theorem \ref{thm:simple} and to ask questions
like those we have asked here, in the case where the two Banach couples
$\left(A_{0},A_{1}\right)$ and $\left(B_{0},B_{1}\right)$ are possibly
different. However for most purposes, and certainly for our purposes
here, there is no loss of generality if we restrict our attention
to the case where $A_{0}=B_{0}$ and $A_{1}=B_{1}$. Let us be a little
more specific about this: Given any any operator $T:A_{0}+A_{1}\to B_{0}+B_{1}$
which satisfies $T(A_{j})\subset B_{j}$ for $j=0,1$, consider the
couple $\left(E_{0},E_{1}\right)=\left(A_{0}\oplus B_{0},A_{1}\oplus B_{1}\right)$
and the operator $S:E_{0}+E_{1}\to E_{0}+E_{1}$ defined by $S(a\oplus b)=\left(0\oplus T(a)\right)$.
Clearly $S\left(E_{j}\right)\subset E_{j}$ for $j=0,1$ and, for
example, $T:\left(A_{0},A_{1}\right)_{\theta,q}\to\left(B_{0},B_{1}\right)_{\theta,q}$
is compact if and only if $S:\left(E_{0},E_{1}\right)_{\theta,q}\to\left(E_{0},E_{1}\right)_{\theta,q}$
is compact.
\end{rem}
\bigskip{}

\textbf{\textit{Acknowledgement.}} We thank Mario Milman for some
very helpful comments.

\section{\label{sec:prelim}Some preliminary results}

\subsection{\label{sec:lambdaf}Some simple nonlinear operators which act on
Banach lattices}

Suppose that $\left(\Omega,\Sigma,\mu\right)$ is an arbitrary measure
space, and that $v:\Omega\to[0,\infty)$ is a fixed measurable function.
In this subsection we will take note of some trivial but useful properties
of three very simple nonlinear operators, which we will denote by
$\Lambda_{v}$, $M_{v}$ and $\widetilde{M}_{v}$. We will define
them by \[
(\Lambda_{v}(f))(\omega)=\min\left\{ \left|f(\omega)\right|,v(\omega)\right\} \,,\]
\[
(M_{v}(f))(\omega)=\max\left\{ \left|f(\omega)\right|,v(\omega)\right\} \,\]
and\begin{equation}
(\widetilde{M}_{v}(f))(\omega)=\max\left\{ \left|f(\omega)\right|,v(\omega)\right\} -v(\omega)\,,\label{eq:defmt}\end{equation}
for all $\omega\in\Omega$ and all measurable functions $f:\Omega\to\mathbb{C}$.

We first claim that each of the three inequalities \begin{equation}
\left|(\Lambda_{v}(f))(\omega)-\Lambda_{v}(g)(\omega)\right|\le\left|f(\omega)-g(\omega)\right|\,,\label{eq:pwmon}\end{equation}

\begin{equation}
\left|(M_{v}(f))(\omega)-M_{v}(g)(\omega)\right|\le\left|f(\omega)-g(\omega)\right|\,\label{eq:mmon}\end{equation}
and\begin{equation}
\left|(\widetilde{M}_{v}(f))(\omega)-\widetilde{M}_{v}(g)(\omega)\right|\le\left|f(\omega)-g(\omega)\right|\,\label{eq:zmon}\end{equation}
holds for all measurable $f:\Omega\to\mathbb{C}$ and $g:\Omega\to\mathbb{C}$
and for all $\omega\in\Omega$.

Of course (\ref{eq:zmon}) is the same as (\ref{eq:mmon}). To prove
(\ref{eq:pwmon}) and (\ref{eq:mmon}) we simply consider the following
four subsets of $\Omega$, namely $\Omega_{-}=\left\{ \omega\in\Omega:\max\{|f(\omega)|,|g(\omega)|\}\le v(\omega)\right\} $,
$\Omega_{+}=\left\{ \omega\in\Omega:\min\{|f(\omega)|,|g(\omega)|\}\ge v(\omega)\right\} $,
$\Omega_{f}=\left\{ \omega\in\Omega:|f(\omega)|\ge v(\omega)\ge|g(\omega)|\right\} $
and $\Omega_{g}=\left\{ \omega\in\Omega:|g(\omega)|\ge v(\omega)\ge|f(\omega)|\right\} $.
Obviously (\ref{eq:pwmon}) and (\ref{eq:mmon}) both hold on each
one of these sets, in each case for some other trivial reasons. Since
$\Omega$ is the union of these sets, the proof of our claim is complete.

Now let $X$ be an arbitrary Banach lattice of (equivalence classes
of) measurable functions on $\left(\Omega,\Sigma,\mu\right)$. Obviously
we have $\left(\Lambda_{v}(f)\right)(\omega)\le\left|f(\omega)\right|$
which implies that \begin{equation}
\Lambda_{v}(X)\subset X\,.\label{eq:flx}\end{equation}
Since $0\le\left(\widetilde{M}_{v}(f)\right)(\omega)=\max\left\{ \left(|f(\omega)|-v(\omega)\right),0\right\} \le\left|f(\omega)\right|$
it also follows that \begin{equation}
\widetilde{M}_{v}(X)\subset X\,.\label{eq:mx}\end{equation}
Furthermore, as an immediate consequence of (\ref{eq:pwmon}) and
(\ref{eq:zmon}), we obtain that $\Lambda_{v}$ and $\widetilde{M}_{v}$
have the Lipschitz norm properties \begin{equation}
\left\Vert \Lambda_{v}(f)-\Lambda_{v}(g)\right\Vert _{X}\le\left\Vert f-g\right\Vert _{X}\,\mbox{and }\left\Vert \widetilde{M}_{v}(f)-\widetilde{M}_{v}(g)\right\Vert _{X}\le\left\Vert f-g\right\Vert _{X}\,\mbox{ for all }f,g\in X\,.\label{eq:lmlp}\end{equation}
Using the fact that $\Lambda_{v}(0)=0$ and $\widetilde{M}_{v}(0)=0$
or the pointwise inequalities mentioned earlier, we also have the
boundedness properties\[
\left\Vert \Lambda_{v}(f)\right\Vert _{X}\le\left\Vert f\right\Vert _{X}\,\mbox{and }\left\Vert \widetilde{M}_{v}(f)\right\Vert _{X}\le\left\Vert f\right\Vert _{X}\,\mbox{ for all }f\in X\,.\]

\subsection{A convenient criterion for showing that an operator is not compact}

The following result will be used for treating most of our examples.
(A slightly different approach will be used for two of the examples
in the appendix.)
\begin{lem}
\label{lem:notcompact}Suppose that the Banach couple $\left(A_{0},A_{1}\right)$
is either $\left(L^{1},L^{\infty}\right)$ or $\left(L^{\infty},L^{1}\right)$
for some arbitrary underlying measure space $\left(\Omega,\Sigma,\mu\right)$.
Suppose that $T$ is a possibly nonlinear map from $A_{0}+A_{1}$
which satisfies $T\left(\left(A_{0},A_{1}\right)_{\theta,q}\right)\subset\left(A_{0},A_{1}\right)_{\theta,q}$
for each $\theta\in(0,1)$ and each $q\in[1,\infty]$.

Suppose that, for each $p\in(1,\infty)$, there exist a sequence $\left\{ E_{N}\right\} _{N\in\mathbb{N}}$
of pairwise disjoint measurable subsets of $\Omega$ and positive
numbers $\nu_{p}$ and $\gamma_{p}$ depending only on $p$, such
that the functions $\psi_{N}=\frac{1}{(\mu(E_{N}))^{1/p}}\chi_{E_{n}}$
satisfy $\gamma_{p}\psi_{N}\le T\left(\psi_{N}\right)\le\psi_{N}$
for each $N>\nu_{p}$. Then, for every $\theta\in(0,1)$ and for every
$q\in[1,\infty]$, the map $T:\left(A_{0},A_{1}\right)_{\theta,q}\to\left(A_{0},A_{1}\right)_{\theta,q}$
is not compact.
\end{lem}
\textit{Proof.} For each choice of $p\in(1,\infty)$, the functions
$\psi_{N}$ defined above (and depending on $p)$ obviously satisfy
$\left\Vert \psi_{N}\right\Vert _{L^{p}}=1$ for all $N$. Furthermore,
whenever $\nu_{p}<N<N'$, we have \[
\left\Vert T(\psi_{N})-T(\psi_{N'})\right\Vert _{L^{p}}=\left\Vert T(\psi_{N})+T(\psi_{N'})\right\Vert _{L^{p}}\ge\gamma_{p}\left\Vert \psi_{N}+\psi_{N'}\right\Vert _{L^{p}}=\gamma_{p}\cdot2^{1/p}\,.\]
This suffices to show that $T$ does not map all bounded subsets of
$L^{p}$ into compact subsets of $L^{p}$.

A slight modification of the preceding argument, using exactly the
same sequence of functions, will now give the corresponding conclusion
for the space $L^{p,q}$, in place of $L^{p}$, for each choice of
$q\in[1,\infty]$. We use the standard quasinorm $\left\Vert f\right\Vert _{L^{p,q}}=\left(\int_{0}^{\infty}\left(t^{1/p}f^{*}(t)\right)^{q}dt/t\right)^{1/q}$
for $L^{p,q}$, with $\left\Vert f\right\Vert _{L^{p,\infty}}=\sup_{t>0}t^{1/p}f^{*}(t)$
when $q=\infty$. We will use two standard properties of non increasing
rearrangements, namely that $(cf)^{*}=cf^{*}$ for each positive constant
$c$ and that $f^{*}\le g^{*}$ whenever $0\le|f|\le|g|$. These,
combined with the fact that the non increasing arrangement of $\chi_{E_{N}}$
is of course the function $\left(\chi_{E_{N}}\right)^{*}(t)=\chi_{[0,\mu(E_{N}))}(t)$,
lead to the following conclusions:

\[
\left\Vert \psi_{N}\right\Vert _{L^{p,q}}=\frac{1}{(\mu(E_{N})^{1/p}}\left(\int_{0}^{\mu(E_{N})}t^{q/p-1}dt\right)^{1/q}=\left(\frac{p}{q}\right)^{1/q}\]
and, whenever $\nu_{p}<N<N'$, \begin{eqnarray*}
\left\Vert T(\psi_{N})-T(\psi_{N'})\right\Vert _{L^{p,q}} & = & \left\Vert T(\psi_{N})+T(\psi_{N'})\right\Vert _{L^{p,q}}\ge\gamma_{p}\left\Vert \psi_{N}+\psi_{N'}\right\Vert _{L^{p,q}}\\
 & \ge & \gamma_{p}\left\Vert \psi_{N}\right\Vert _{L^{p,q}}=\gamma_{p}\cdot\left(\frac{p}{q}\right)^{1/q}\,.\end{eqnarray*}
When $q=\infty$ we obtain the same conclusions, with $\left(\frac{p}{q}\right)^{1/q}$
replaced by $1$.

Thus, in all cases, the sequence $\left\{ T(\psi_{N})\right\} _{N\in\mathbb{N}}$
cannot have a subsequence which converges in $L^{p,q}$.

To complete the proof of the lemma it remains to recall that, for
our choices of the couple $\left(A_{0},A_{1}\right)$, and for each
$\theta\in(0,1)$ and each $q\in[1,\infty]$, the space $(A_{0},A_{1})_{\theta,q}$
always coincides with $L^{p,q}$ for some $p\in(1,\infty)$, to within
equivalence of quasinorms. In fact (see e.g., Theorem 5.3.1 on p.~113
of \cite{bl}) $(L^{1},L^{\infty})_{\theta,q}$ and $\left(L^{\infty},L^{1}\right)_{\theta,q}$
coincide respectively with $L^{\frac{1}{1-\theta},q}$ and $L^{\frac{1}{\theta},q}$.
$\qed$

\section{\label{sec:onesided}A one-sided compactness assumption on the bigger
space is not sufficient}

In this section we shall present our first counterexample, a rather
simple nonlinear operator $T_{1}:L^{1}\to L^{1}$ which has the properties
{[}i{]}, {[}ii{]}, {[}iii$_{1}]$ and {[}iv{]}, for $A_{0}=L^{\infty}$
and $A_{1}=L^{1}$ on the measure space $[0,1]$.

For each $n\in\mathbb{N}$ let $I_{n}$ be the open interval $(2^{-n},2^{-n+1})$.
Define the function $v:[0,1]\to[0,\infty)$ by \[
v=\sum_{n=1}^{\infty}\frac{2^{n}}{n^{2}}\chi_{I_{n}}\]
and let $Q:L^{1}\to L^{1}$ be the conditional expectation operator
defined by \[
Qf=\sum_{n=1}^{\infty}\frac{1}{\left|I_{n}\right|}\int_{I_{n}}f(x)dx\cdot\chi_{I_{n}}\,.\]
Our operator $T_{1}$ is given by the formula \[
T_{1}(f)=\min\left\{ \left|Qf\right|,v\right\} \,\mbox{for all }f\in L^{1}.\]
In other words, $T_{1}$ is the composition of operators $T_{1}=\Lambda_{v}\circ Q$.
It obviously satisfies property {[}i{]}. Since both $Q$ and $\Lambda_{v}$
both have property {[}ii{]}, (cf.~(\ref{eq:flx}) and (\ref{eq:lmlp}))
so does their composition $T_{1}$.

Let $H$ be the set of all functions $f:[0,1]\to[0,\infty)$ of the
form $f=\sum_{n=1}^{\infty}\alpha_{n}\chi_{I_{n}}$ where each of
the constants $\alpha_{n}$ satisfies $0\le\alpha_{n}\le\frac{2^{n}}{n^{2}}$.
The convergence of the series $\sum_{n=1}^{\infty}\frac{1}{n^{2}}$
ensures that, for each $\epsilon>0$, there exists $N_{\epsilon}$
such that $\sum_{n\ge N_{\epsilon}}\frac{2^{n}}{n^{2}}\left|I_{n}\right|<\epsilon$.
It follows readily that $H$ is a compact subset of $L^{1}$. Since
$T_{1}(L^{1})\subset H$, we see that $T_{1}$ certainly has property
{[}iii$_{1}${]}.

Finally, we have to show that $T_{1}$ has property {[}iv{]}. Let
us choose an arbitrary number $p\in(1,\infty)$. For each $N\in\mathbb{N}$,
let $\psi_{N}=2^{N/p}\chi_{I_{N}}$. There exists an integer $\sigma_{p}$
depending on $p$ such that, for all $N\ge\sigma_{p}$, we have $2^{N/p}\le\frac{2^{N}}{N^{2}}$
and therefore $T(\psi_{N})=\psi_{N}$. Clearly we can now apply Lemma
\ref{lem:notcompact}, with $E_{N}=I_{N}$ and $\gamma_{p}=1$ and
$\nu_{p}=\sigma_{p}$, to obtain property {[}iv{]}.

\section{\label{sec:smallonesided}A one-sided compactness assumption on the
smaller space is not sufficient}

In this section we present our second counterexample. It uses the
same couple $\left(A_{0},A_{1}\right)=\left(L^{\infty},L^{1}\right)$
with the same underlying measure space $[0,1]$ and same sequence
of intervals $\left\{ I_{n}\right\} _{n\in\mathbb{N}}$ and the same
conditional expectation operator $Q$ as the example of the previous
section. This time, instead of property {[}iii$_{1}${]} of the above
list, we will obtain property {[}iii$_{0}${]}, together with {[}i{]},
{[}ii{]} and {[}iv{]}.

Let $w:[0,1]\to[0,\infty)$ be the function \[
w=\sum_{n=1}^{\infty}n\chi_{I_{n}}\]
 and let $T_{2}$ be the nonlinear operator $T=\widetilde{M}_{w}\circ Q$.
I.e., we set \[
T_{2}(f)=\sum_{n=1}^{\infty}\left(\max\left\{ n,\frac{1}{\left|I_{n}\right|}\left|\int_{I_{n}}f(x)dx\right|\right\} -n\right)\chi_{I_{n}}\,\mbox{for each }f\in L^{1}\,.\]
Both of the operators $\widetilde{M}_{w}$ and $Q$ satisfy properties
{[}i{]} and {[}ii{]} (cf.~(\ref{eq:mx}) and (\ref{eq:lmlp})). Therefore,
so does their composition $T_{2}$.

Now, to establish {[}iii$_{0}${]}, let $A$ be an arbitrary bounded
subset of $L^{\infty}$. Choose an integer $N$ such that $\left\Vert f\right\Vert _{L^{\infty}}\le N$
for all $f\in A$. Then, for each $f\in A$, the function $T_{2}(f)$
vanishes at every point of the set $\bigcup_{n\ge N}I_{n}$ and is
constant and bounded by $N-n$ on each of the intervals $I_{n}$ for
$1\le n<N$. Thus $T_{2}(A)$ is contained in the set \[
H_{N}=\left\{ f=\sum_{k=1}^{N-1}\gamma_{n}\chi_{I_{n}}:0\le\gamma_{n}\le N\right\} \]
which is of course compact in $L^{\infty}$.

Finally we will show that property {[}iv{]} holds. For each choice
of $p\in(1,\infty)$ we will use the same functions $\psi_{N}=2^{N/p}\chi_{I_{N}}$
as we used in Section \ref{sec:onesided}. This time we have $T_{2}(\psi_{N})=\left(\max\left\{ N,2^{N/p}\right\} -N\right)\chi_{I_{N}}=\max\left\{ 2^{N/p}-N,0\right\} \chi_{I_{N}}$
for each $N\in\mathbb{N}$. There exists a positive integer $\tau_{p}$
such that, whenever $N>\tau_{p}$, we have $N\le2^{(N-1)/p}$ and
therefore also \begin{equation}
2^{N/p}\left(1-2^{-1/p}\right)\chi_{I_{N}}\le(2^{N/p}-N)\chi_{I_{N}}=T_{2}(\psi_{N})\,.\label{eq:ftt}\end{equation}

These properties enable us to obtain property {[}iv{]} by applying
Lemma \ref{lem:notcompact}, with $E_{N}=I_{N}$ as before, but this
time with $\gamma_{p}=1-2^{-1/p}$ and $\nu_{p}=\tau_{p}$.

\section{\label{sec:twosided}Even a two sided compactness assumption is not
sufficient}

We can now combine the operators $T_{1}$ and $T_{2}$ of the preceding
two examples to obtain our main counterexample which has all of the
properties {[}i{]}, {[}ii{]}, {[}iii$_{0}${]}, {[}iii$_{1}${]} and
{[}iv{]}. Our new operator will simply be their composition $T_{3}=T_{2}\circ T_{1}$.

The spaces will be, as before, $A_{0}=B_{0}=L^{\infty}$ and $A_{1}=B_{1}=L^{1}$.
The intervals $I_{n}$ will also be as before.

Since both $T_{1}$ and $T_{2}$ satisfy properties {[}i{]} and {[}ii{]},
so does their composition.

In Section \ref{sec:onesided} we saw that $T_{1}(L^{1})$ is contained
in the compact subset $H$ of $L^{1}$. Since $T_{2}$ is a Lipschitz
and therefore continuous map of $L^{1}$ into itself, it must map
$H$ into another compact subset of $L^{1}$. Thus $T_{3}$ has property
{[}iii$_{1}${]}.

Since $T_{1}$ has properties {[}i{]} and {[}ii{]}, the set $T_{1}(A)$
is bounded in $L^{\infty}$ whenever $A$ is a bounded subset of $L^{\infty}$.
Since $T_{2}$ has property {[}iii$_{0}${]}, the set $T_{3}(A)=T_{2}\left(T_{1}(A)\right)$
must be relatively compact in $L^{\infty}$, and this shows that $T_{3}$
has property {[}iii$_{0}${]}.

Now, for property {[}iv{]}, we choose an arbitrary $p\in(1,\infty)$
and once more consider the functions $\psi_{N}=2^{N/p}\chi_{I_{N}}$.
If $N>\max\left\{ \sigma_{p},\tau_{p}\right\} $ we have $T_{1}(\psi_{N})=\psi_{N}$
and so $T_{3}(\psi_{N})=T_{2}(\psi_{N})$. Then (\ref{eq:ftt}) gives
us that \[
2^{N/p}\left(1-2^{-1/p}\right)\chi_{I_{N}}\le T_{2}(\psi_{N})=T_{3}(\psi_{N})\,.\]
This enables us to apply Lemma \ref{lem:notcompact} one more time,
this time with $E_{N}=I_{N}$ and $\gamma_{p}=1-2^{-1/p}$, as in
Section \ref{sec:smallonesided}, but now with $\nu_{p}=\max\left\{ \sigma_{p},\tau_{p}\right\} $.
This gives us that $T_{3}$ has property {[}iv{]}.

\section{\label{sec:appendix}Appendix - Some additional counterexamples}

\subsection{\label{sub:sequence}An example for a couple of sequence spaces }

In this subsection we will describe another simple operator $T_{4}$
which has the same properties {[}i{]}, {[}ii{]}, {[}iii$_{1}${]}
and {[}iv{]} as the operator $T_{1}$ of Section \ref{sec:onesided}.
But here, in contrast to all the other examples in this paper, the
couple $\left(A_{0},A_{1}\right)$ will be the couple of sequence
spaces $\left(\ell^{1},\ell^{\infty}\right)$. Note that here, as
in all the other examples and as in Theorem \ref{thm:simple}, we
still have $A_{0}$ continuously embedded in $A_{1}$.

Let $v:[0,\infty)\to[0,\infty)$ be the function $v=\sum_{n=1}^{\infty}\frac{1}{n}\chi_{[2^{n},2^{n+1})}$.
In fact we will use the restriction of $v$ to the set $\mathbb{N}$,
i.e., the sequence $\left\{ v(n)\right\} _{n\in\mathbb{N}}$.

Let $T_{4}:\ell^{\infty}\to\ell^{\infty}$ be the operator which maps
each bounded sequence $\alpha=\left\{ \alpha_{n}\right\} _{n\in\mathbb{N}}$
to the sequence $T(\alpha)$ which is defined by the formula \[
\left(T_{4}(\alpha)\right)_{n}=\min\left\{ \left|\alpha_{n}\right|,v(n)\right\} \,.\]

In other words, we have chosen $T_{4}$ to be the operator $\Lambda_{v}$
of Subsection \ref{sec:lambdaf}, in the case where the underlying
measure space is $\mathbb{N}$ equipped with counting measure, and
the function $v:\mathbb{N}\to[0,\infty]$ is given by $v(n)$ defined
as above. Property {[}i{]} is immediate. By (\ref{eq:flx}) and (\ref{eq:lmlp}),
we also immediately have that $T_{4}$ has property {[}ii{]}.

Since $\lim_{n\to\infty}v(n)=0$, the set $H=\left\{ \alpha\in\ell^{\infty}:\left|\alpha_{n}\right|\le v(n)\,\mbox{for all }n\in\mathbb{N}\right\} $
is a compact subset of $c_{0}$ and therefore also of $\ell^{\infty}$.
Since $T_{4}$ maps $\ell^{\infty}$ onto $H$ we certainly have that
$T_{4}$ maps every bounded subset of $\ell^{\infty}$ into a relatively
compact subset of $\ell^{\infty}$. This establishes property {[}iii$_{1}]$.

As previously, we will use Lemma \ref{lem:notcompact} to establish
property {[}iv{]}. This time we choose our sequence $\left\{ E_{N}\right\} _{N\in\mathbb{N}}$
of pairwise disjoint subsets of the underlying measure space by setting
$E_{N}=\left\{ n\in\mathbb{N}:2^{N}\le n<2^{N+1}\right\} $. Then,
after fixing $p\in(1,\infty)$, since $E_{N}$ contains $2^{N}$ points,
we have to choose $\psi_{N}=\left\{ \left(\psi_{N}\right)_{n}\right\} _{n\in\mathbb{N}}$
to be the sequence defined by \[
\left(\psi_{N}\right)_{n}=\left\{ \begin{array}{ccc}
2^{-N/p} & , & n\in E_{N}\\
0 & , & n\in\mathbb{N}\backslash E_{N}\end{array}\right.\,.\]
There exists an integer $\nu_{p}$ depending only on $p$ with the
property that, each integer $N>\nu_{p}$ satisfies $1/N>2^{-N/p}$
and therefore also $T_{4}\left(\psi_{N}\right)=\psi_{N}$. So we can
now apply Lemma \ref{lem:notcompact} for this particular choice of
$\nu_{p}$ and for $\gamma_{p}=1$.

\subsection{\label{sub:elaborate}A more elaborate example for the couple $\left(A_{0},A_{1}\right)=\left(B_{0},B_{1}\right)=(L^{\infty},L^{1})$}

Our example in this subsection is considerably more complicated than
all preceding counterexamples. As in our first three counterexamples,
we will again take $\left(A_{0},A_{1}\right)=\left(B_{0},B_{1}\right)=\left(L^{\infty},L^{1}\right)$
where the underlying measure space is the unit interval $[0,1]$ equipped
with Lebesgue measure. Our operator $T_{5}$ will satisfy properties
{[}i{]}, {[}ii{]} and {[}iii$_{0}${]}. But, instead of obtaining
property {[}iv{]}, we will only show that $T_{5}$ is not compact
on the space $\left(A_{0},A_{1}\right)_{1/p,p}=\left(B_{0},B_{1}\right)_{1/p,p}=L^{p}$,
for one particular choice of $p\in(1,\infty)$. Our construction of
$T_{5}$ with these properties will work for any value that we please
of $p$ in $\left(1,\infty\right)$. But it depends on our choice
of that $p$.

Once we have made our choice of $p$, we fix a sequence of numbers
$\left\{ m_{n}\right\} _{n\ge0}$ defined by $m_{0}=0$ and $m_{n}=\sum_{k=1}^{n}2^{k(p+1)}$
for $n\ge1$. We also fix a sequence $\left\{ I_{N}\right\} _{N\in\mathbb{N}}$
of pairwise disjoint open subintervals of $[0,1]$ such that the length
of $I_{N}$ is $2^{-Np}$ for each $N$. Note that $\sum_{N=1}^{\infty}2^{-Np}=\frac{2^{-p}}{1-2^{-p}}=\frac{1}{2^{p}-1}<1$
so the interval $[0,1]$ is sufficiently large to accomodate such
a sequence. Our operator $T_{5}$ will be defined as a pointwise supremum
of a sequence of functions, by the formula \begin{equation}
T_{5}(f)=\sup_{N\in\mathbb{N}}S_{N}\left(\frac{1}{2^{-Np}}\int_{I_{N}}\left|f(x)\right|dx\cdot\chi_{I_{N}}\right)\,\mbox{for each }f\in L^{1}[0,1]\mbox{\ensuremath{\,}, }\label{eq:deft}\end{equation}
where, for each $N$, we take $S_{N}$ to be an appropriately defined
nonlinear operator acting on the one dimensional space of functions
$\left\{ c\chi_{I_{N}}:c\in\mathbb{C}\right\} $.

\subsubsection{\label{sub:auxiliary}Construction of the auxiliary operators $S_{N}$.}

In this subsubsection we carry out the major step of constructing
each of the operators $S_{N}$ and then obtain some of their properties
which will be needed later to show that $T_{5}$ has all the required
properties. We will proceed somewhat indirectly. We first choose some
arbitrary but fixed positive integer $N$. Since we will have other
subscripts and superscripts in our construction, let us suppress mention
of $N$ for the moment, and simplify the notation by writing $w$
for the length (or width) $2^{-Np}$ of the interval $I_{N}$.

We introduce the numerical sequence $\left\{ h_{n}\right\} _{n\ge0}$,
defined by \[
h_{n}=\sqrt{w^{2}+\frac{w}{2^{p(n+1)}}}-w\mbox{ \ensuremath{\,}.}\]
Note that $\left\{ h_{n}\right\} _{n\ge0}$ is a strictly positive
and strictly decreasing sequence. Next we define two more numerical
sequences $\left\{ y_{n}\right\} _{n\ge1}$ and $\left\{ \lambda_{n}\right\} _{n\ge0}$
by setting $\lambda_{0}=0$ and, for each $n\ge1$, setting \[
y_{n}=\frac{2^{n-1}}{1+\frac{h_{n-1}}{2w}}\mbox{ and }\lambda_{n}=\sum_{k=1}^{n}y_{k}\,.\]

The properties of $\left\{ h_{n}\right\} _{n\ge0}$ ensure that $0<y_{1}\le y_{n}\le y_{n+1}=\lambda_{n+1}-\lambda_{n}$
for each $n\in\mathbb{N}$. Therefore $\lim_{n\to\infty}\lambda_{n}=\infty$
and we can express the interval $(0,\infty)$ as a union of pairwise
disjoint intervals \begin{equation}
(0,\infty)=\bigcup_{n\in\mathbb{N}}(\lambda_{n-1},\lambda_{n}]\,.\label{eq:lnz}\end{equation}
This also means that, for each $t>0$, there exists a unique positive
integer $\nu(t)$ such that \begin{equation}
\lambda_{\nu(t)-1}<t\le\lambda_{\nu(t)}\,.\label{eq:defnu}\end{equation}
We also want to define $\nu(t)$ when $t=0$. We can take $\nu(0)=0$.
(Then we can also arrange to have (\ref{eq:defnu}) also hold when
$t=0$, provided we define $\lambda_{n}$ when $n=-1$ and choose
$\lambda_{-1}$ to be some negative number.)

We are now going to construct a family $\left\{ E(t)\right\} _{t\ge0}$
of subsets of $\mathbb{R}^{2}$. At first we will describe the set
$E(t)$ only for those numbers $t$ which coincide with some element
of the sequence $\left\{ \lambda_{n}\right\} _{n\ge0}$. We will use
the abbreviated notation $E_{n}=E\left(\lambda_{n}\right)$ for these
particular sets.

In each case where $n\ge1$ the set $E_{n}$ is the union of a (solid)
closed rectangle $R_{n}$ whose sides are parallel to the axes, with
a (solid) closed triangle $\Delta_{n}$ located on the right side
of the rectangle. The vertices of $R_{n}$ are $\left(0,\lambda_{n-1}\right)$,
$\left(0,\lambda_{n}\right)$, $\left(w,\lambda_{n-1}\right)$ and
$\left(w,\lambda_{n}\right)$. These last two points are also vertices
of $\Delta_{n}$ and the third vertex of $\Delta_{n}$ is the point
$(w+h_{n-1},\lambda_{n-1}-m_{n-1}h_{n-1})$.

The following very approximate picture of the set $E_{n}$ (for some
$n\ge2$) may be helpful.

\begin{center}
\includegraphics[scale=0.3]{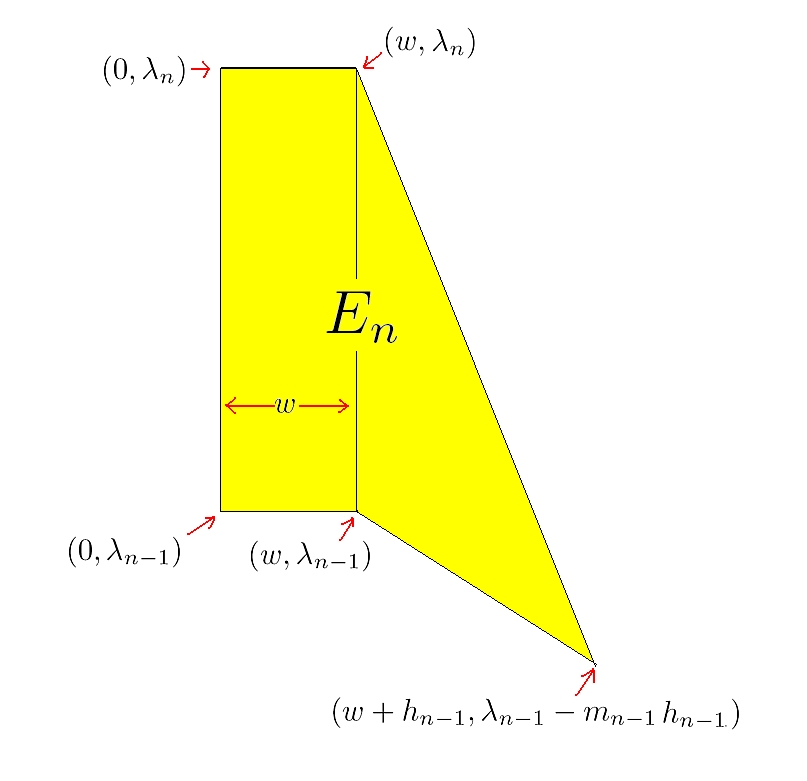}
\par\end{center}

The formulae for the various preceding sequences which are used to
define these vertices of $E_{n}$ are not quite as mysterious as they
may first appear to be. Their choices have been completely determined
by the need to ensure that the area of $E_{n}$ and the slopes of
two non vertical sides of $\partial\Delta_{n}$ are given by some
rather simple formulae, which we shall now obtain.

We first determine the slopes of the two non vertical sides of $\partial\Delta_{n}$.
The slope of the lower one of these sides is of course $-m_{n-1}$.
The slope of the upper side equals \begin{eqnarray}
\frac{(\lambda_{n-1}-m_{n-1}h_{n-1})-\lambda_{n}}{h_{n-1}} & = & -\frac{y_{n}+m_{n-1}h_{n-1}}{h_{n-1}}\nonumber \\
 & = & -\frac{y_{n}}{h_{n-1}}-m_{n-1}\,.\label{eq:ityc}\end{eqnarray}
Now \begin{eqnarray*}
\frac{y_{n}}{h_{n-1}} & = & \frac{2^{n-1}}{h_{n-1}+\frac{h_{n-1}^{2}}{2w}}=\frac{2^{n}w}{2wh_{n-1}+h_{n-1}^{2}}\\
 & = & \frac{2^{n}w}{h_{n-1}\left(h_{n-1}+2w\right)}=\frac{2^{n}w}{\left(\sqrt{w^{2}+\frac{w}{2^{pn}}}-w\mbox{ }\right)\left(\sqrt{w^{2}+\frac{w}{2^{pn}}}+w\mbox{ }\right)}\\
 & = & \frac{2^{n}w}{w^{2}+\frac{w}{2^{pn}}-w^{2}}=2^{n+pn}=2^{n(p+1)}=m_{n}-m_{n-1}.\end{eqnarray*}
Substituting this in (\ref{eq:ityc}), we see that the slope of the
upper side equals $-m_{n}$.

We will use the usual notation $\left|E\right|$ for the area or two
dimensional Lebesgue measure of any given measurable subset $E$ of
$\mathbb{R}^{2}$. In particular, the area of $E_{n}$ is given by
\[
\left|E_{n}\right|=y_{n}w+\frac{1}{2}y_{n}h_{n-1}=\frac{2^{n-1}}{1+\frac{h_{n-1}}{2w}}\left(w+\frac{h_{n-1}}{2}\right)=2^{n-1}w\,.\]

Here is another very approximate picture, this time of the sets $E_{1}$
and $E_{2}$.

\begin{center}
\includegraphics[scale=0.25]{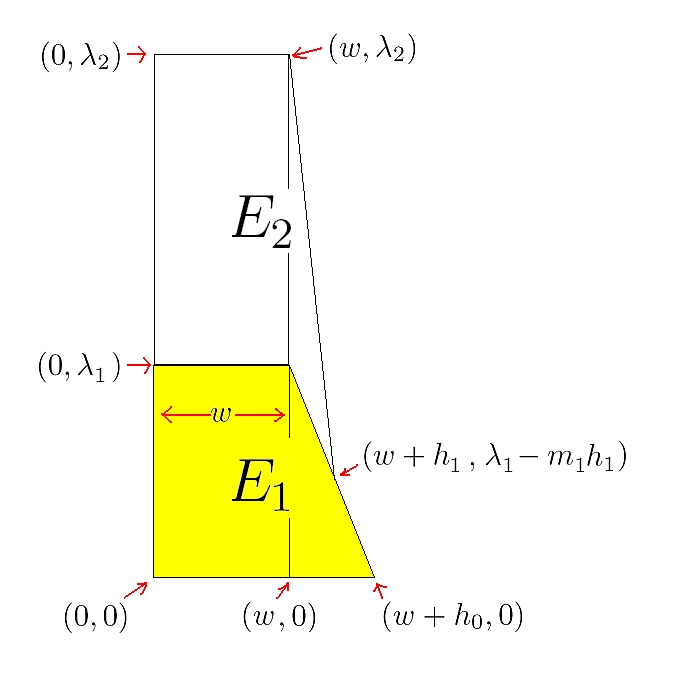}
\par\end{center}

Since $m_{0}=0$ and $\lambda_{0}=0$, we obtain that the set $E_{1}$
is a trapezium (in British terminology) or a trapezoid (in American
terminology) whose base is the line segment from $(0,0)$ to $(w+h_{0},0)$
and which lies entirely in the closed upper half plane.

For each $n\in\mathbb{N}$ the set $E_{n+1}$ fits exactly on top
of the set $E_{n}$ with no overlap. This is indicated by the above
picture when $n=1$ and by the following (approximate) pictore for
$n\ge2$.

\begin{center}
\includegraphics[scale=0.25]{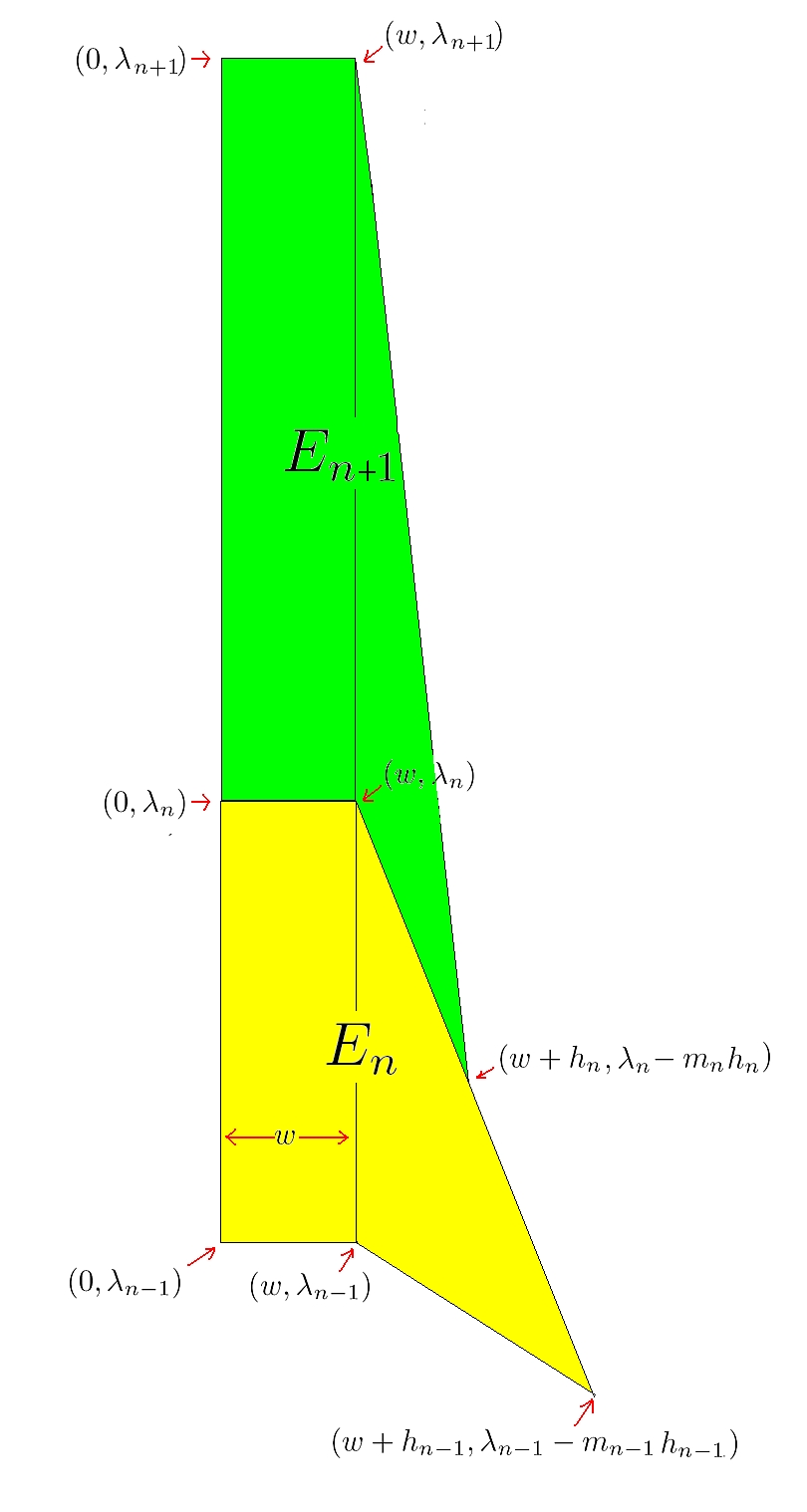}
\par\end{center}

To state this more precisely, we first note that (obviously) the upper
horizontal part of $\partial R_{n}$, coincides with the lower horizontal
part of $\partial R_{n+1}$. Then the upper non vertical side of $\partial\Delta_{n}$
and the lower non vertical side of $\partial\Delta_{n+1}$ both have
the same slope $-m_{n}$ and the same left endpoint $(w,\lambda_{n})$.
Since $0<h_{n}<h_{n-1}$, we see that the first of these sides strictly
contains the second.

We still have to define the set $E_{n}$ for the case where $n=0$.
We will let $E_{0}$ be the non negative $x$ axis, i.e., $E_{0}=\left\{ (x,0):x\ge0\right\} $.

Now we can extend our definition of $E(\lambda_{n})=E_{n}$ to define
the sets $E(t)$ also for those $t\ge0$ which do not coincide with
any $\lambda_{n}$. In view of (\ref{eq:lnz}), this means we have
to define $E(t)$ for each $t$ in the interval $(\lambda_{n-1},\lambda_{n})$
and to do this for each $n\in\mathbb{N}$.

So let us fix some arbitrary $n\in\mathbb{N}$ and consider all numbers
$t\in(\lambda_{n-1},\lambda_{n})$. Note that all these numbers satisfy
$\nu(t)=n$, where $\nu(t)$ is the integer defined in (\ref{eq:defnu})
above. For each $t$ in this interval, the set $E(t)$ is the subset
of $E_{\nu(t)}=E_{n}$ shown (approximately) as the shaded area in
the following picture.

\begin{center}
\includegraphics[scale=0.25]{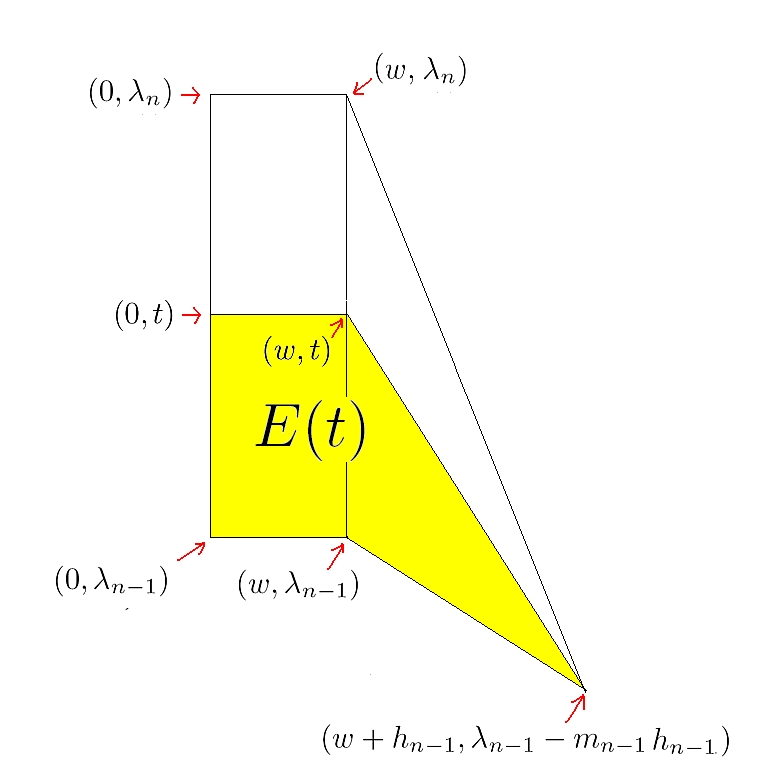}
\par\end{center}

More precisely, $E(t)$ consists of all those points of $E_{n}$ which
lie on or below two particular straight lines, the horizontal line
$y=t$ and the line which passes through the points $(w,t)$ and $(w+h_{n-1},\lambda_{n-1}-m_{n-1}h_{n-1})$.
In other words, $E(t)$ is defined exactly like $E_{n}$, except that
the two uppermost vertices $(0,\lambda_{n})$ and $(w,\lambda_{n})$
are replaced by the two lowered points $(0,t)$ and $(w,t)$. For
later purposes we note that the slope $\sigma(t)$ of the oblique
line which forms part of the upper boundary of $E(t)$ is negative
and its value lies between the values of the slopes of the two non
vertical sides of the triangle $\partial\Delta_{n}$. Thus, from our
previous calculations of slopes, we have that $m_{n-1}<\left|\sigma(t)\right|<m_{n}$.
It will be convenient to rewrite this as a formula which will be valid
for all $t>0$, namely \begin{equation}
m_{\nu(t)-1}<\left|\sigma(t)\right|<m_{\nu(t)}\,.\label{eq:msig}\end{equation}

The area of $E(t)$ is given by the formula \begin{equation}
\left|E(t)\right|=\left(w+\frac{1}{2}h_{n-1}\right)\left(t-\lambda_{n-1}\right)\,\mbox{for each }t\in(\lambda_{n-1},\lambda_{n})\,.\label{eq:rog}\end{equation}
Or, in other words, $\left|E(t)\right|=\left(w+\frac{1}{2}h_{\nu(t)-1}\right)\left(t-\lambda_{\nu(t)-1}\right)$
for each $t\ge0$ which is not an element of the sequence $\left\{ \lambda_{n}\right\} _{n\ge0}$.
Thus we see that, for each $n\in\mathbb{N}$, the function $t\mapsto\left|E(t)\right|$
is a positive strictly increasing affine function on the open interval
$\left(\lambda_{n-1},\lambda_{n}\right)$ and its limits at $\lambda_{n-1}$
and $\lambda_{n}$ (one sided limits with respect to this interval)
are $0$ and $\left|E_{n}\right|=2^{n-1}w$ respectively.

Our next step is to use the family of sets $\left\{ E(t)\right\} _{t\ge0}$
to define another family of planar sets which we will denote by $\left\{ G(t)\right\} _{t\ge0}$.

Analogously to our handling of the family $\left\{ E(t)\right\} _{t\ge0}$,
we shall begin by defining the sets $G(t)$ when $t=\lambda_{n}$
for some integer $n$ and by using the notation $G_{n}=G(\lambda_{n})$.
For each integer $n\ge0$ we let $G(\lambda_{n})=G_{n}=\bigcup_{k=0}^{n}E_{k}$.
Then $\left|G_{n}\right|$ is of course equal to the sum of the areas
of (the interiors of) the non overlapping sets $E_{k}$ and thus it
is given by $\sum_{k=1}^{n}2^{k-1}w=(2^{n}-1)w$. Since $G_{n}$ contains
a rectangle of width $w$ and height $\lambda_{n}$ we clearly have
\begin{equation}
\lambda_{n}\le2^{n}-1.\label{eq:efl}\end{equation}
Note that the formula $(2^{n}-1)w$ for $\left|G_{n}\right|$ and
the estimate (\ref{eq:efl}) for $\lambda_{n}$ both hold also in
the trivial case where $n=0$.

For the remaining values of $t\ge0$, i.e., those which do not coincide
with any $\lambda_{n}$, we set \[
G(t)=G(\lambda_{\nu(t)-1})\cup E(t)\]
where the integer $\nu(t)$ is defined as before. In other words,
we have $G(t)=\left(\bigcup_{k=0}^{\nu(t)-1}E_{k}\right)\cup E(t)$.
It is clear that \begin{equation}
G(t)\subset G(t')\,\mbox{whenever }0\le t\le t'\,.\label{eq:gmon}\end{equation}
It is also clear that $t\mapsto\left|G(t)\right|$ is a continuous
strictly increasing and in fact piecewise affine function on $[0,\infty)$
which satisfies \begin{equation}
\left|G\left(\lambda_{n}\right)\right|=\left(2^{n}-1\right)w\ \mbox{for each integer }n\ge0\,.\label{eq:edm}\end{equation}
In particular this gives us $\left|G(0)\right|=0$ and we also have
$\lim_{t\to\infty}\left|G(t)\right|=\infty$. All these properties
guarantee the existence of an inverse function, namely a continuous
strictly increasing and in fact piecewise affine function $\gamma:[0,\infty)\to[0,\infty)$
which has the property \begin{equation}
\left|G\left(\gamma(s)\right)\right|=s\mbox{ \mbox{for each }}s\ge0\label{eq:fgg}\end{equation}
which can of course also be equivalently expressed as \begin{equation}
\gamma\left(\left|G(t)\right|\right)=t\mbox{ for each }t\ge0\,.\label{eq:sgg}\end{equation}
We will need two more more special properties of $\gamma$. In particular
we remark, using (\ref{eq:edm}) and (\ref{eq:sgg}), that \begin{equation}
\gamma\left((2^{n}-1)w\right)=\lambda_{n}\,\mbox{for each integer }n\ge0\,.\label{eq:twsu}\end{equation}
We also remark that, by (\ref{eq:rog}), each line segment of the
graph of the function $t\mapsto\left|G(t)\right|$ has a positive
slope which is strictly greater than $w$. This means that each line
segment of the graph of the inverse function $s\mapsto\gamma(s)$
has a positive slope which strictly is less than $1/w$. This in turn
ensures that $\gamma$ satisfies the Lipschitz condition \begin{equation}
\left|\gamma(s)-\gamma(s')\right|\le\frac{1}{w}\left|s-s'\right|\,\mbox{for all }s,s'\in[0,\infty)\,.\label{eq:gamlip}\end{equation}
\medskip

We are now ready to define a special function of two variables $g:[0,\infty)\times[0,\infty)\to[0,\infty)$
by the formula \[
g(x,t)=\sup\left\{ y:(x,y)\in G(t)\right\} \,.\]
In other words, for each fixed $t\ge0$, we take $x\mapsto g(x,t)$
to be the function of one variable whose graph is the upper edge of
the set $G(t)$. Since, for all $t>0$, the sets $E(t)$ are all contained
in the strip $\left\{ (x,y):0\le x\le w+h_{0}\right\} $, and since
$G(0)=E_{0}$ is simply the non negative $x$-axis, we see that $g(x,t)=0$
for all $x>w+h_{0}$ and all $t>0$, and also that $g(x,0)=0$ for
all $x\ge0$.

For each fixed $t>0$ we can equivalently reformulate the definition
of the function $x\mapsto g(x,t)$ by declaring it to be the continuous
piecewise affine function which vanishes on the interval $(w+h_{0},\infty)$
and which has a constant derivative on each of the $\nu(t)+1$ intervals
$[0,w)$, $(w,w+h_{\nu(t)-1})$, $\left(w+h_{\nu(t)-1},w+h_{\nu(t)-2}\right),$
... $(w+h_{1},w+h_{0})$ and whose values at the end points of these
intervals are $g(0,t)=t$, $g(w,t)=t$ and $g(w+h_{k},t)=\lambda_{k}-m_{k}h_{k}$
for $k=\nu(t)-1,\nu(t)-2,.....,0$. Note that this formulation is
valid whether or not $t$ is one of the numbers $\lambda_{n}$.

Obviously the derivative $\frac{\partial g}{\partial x}(x,t)$ is
zero for all $x$ in the first interval $[0,w)$. For all $x$ in
the second interval $(w,w+h_{\nu(t)-1})$ it is clear from our preceding
remarks and calculations that $\frac{\partial g}{\partial x}(x,t)$
equals either $\sigma(t)$ or $-m_{\nu(t)}$, depending on whether
$t<\lambda_{\nu(t)}$ or $t=\lambda_{\nu(t)}$. The values of this
derivative on the remaining intervals of the list are, respectively,
$-m_{\nu(t)-1}$,...., $-m_{1}$. In view of (\ref{eq:msig}) and
the fact that $0\le m_{n-1}<m_{n}$ for each $n\in\mathbb{N}$, we
deduce that $\left|\frac{\partial g}{\partial x}(x,t)\right|\le m_{\nu(t)}$
for each $t\ge0$ and for each $x\ge0$ which does not coincide with
any of the {}``cusp'' points $w$ and $w+h_{k}$, $k=\nu(t)-1,\nu(t)-2,.....,0$.
This means that $g$ satisfies the Lipschitz condition \begin{equation}
\left|g(x,t)-g(x',t)\right|\le m_{\nu(t)}\left|x-x'\right|\mbox{ for all non negative }x,\, x,\,\mbox{ and }t\,.\label{eq:lipc}\end{equation}

It is clear that the integral of $g$ for each fixed $t$ has to satisfy
\begin{equation}
\int_{0}^{\infty}g(x,t)dx=\int_{0}^{w+h_{0}}g(x,t)dx=\left|G(t)\right|\,.\label{eq:mbet}\end{equation}
 This means that $\int_{0}^{w}g(x,t)dx\le\left|G(t)\right|$, and
since $g(x,t)=g(0,t)=t=\sup_{s\ge0}g(s,t)$ for each $x\in[0,w]$,
we deduce that \begin{equation}
\sup_{x\ge0}g(x,t)=t\le\frac{\left|G(t)\right|}{w}\,.\label{eq:verp}\end{equation}
Since this tells us that $tw\le\left|G(t)\right|$ we can apply the
monotonicity of $\gamma$ and (\ref{eq:sgg}) to obtain that \begin{equation}
\gamma\left(tw\right)\le\gamma\left(\left|G(t)\right|\right)=t\ \mbox{for each }t\ge0\,.\label{eq:mma}\end{equation}

We also need some facts about the function $g$ considered as a function
of $t$ for fixed values of $x$. First it is clear from (\ref{eq:gmon})
that $t\mapsto g(x,t)$ is a non decreasing function for each fixed
$x$. Then we want to show that the function $g$ satisfies a second
kind of Lipschitz condition. We claim that\begin{equation}
\left|g(x,t)-g(x,t')\right|\le\left|t-t'\right|\ \mbox{for all non negative }x,\, t\,\mbox{and }t'\,.\label{eq:slipc}\end{equation}

We may of course suppose without loss of generality that $0\le t<t'$,
and then, in view of the monotonicity of $t\mapsto g(x,t)$, the condition
(\ref{eq:slipc}) is the same as \begin{equation}
0\le g(x,t')-g(x,t)\le t'-t\ \mbox{for all }x\ge0\mbox{ and }0\le t<t'\,.\label{eq:xlip}\end{equation}

Our first step will be to prove (\ref{eq:xlip}) in the special case
where $t$ and $t'$ are both numbers in the same interval $[\lambda_{n-1},\lambda_{n}]$.
It is clear from the definitions of $g$ and of the sets $E(t)$ and
$G(t)$, that $g(x,t')-g(x,t)=t'-t$ for all $x\in[0,w]$ and that
$g(x,t')-g(x,t)<t'-t$ for all $x\in(w,w+h_{n-1})$. We also have
$g(x,t)=g(x,t')$ for all $x\ge w+h_{n-1}$. Together, these three
properties give us (\ref{eq:xlip}) in this case.

Our second and final step will be to show that (\ref{eq:xlip}) in
fact holds for all $0\le t<t'$ in the remaining case where $t$ and
$t'$ are not in the same interval $[\lambda_{n-1},\lambda_{n}]$
for any $n\in\mathbb{N}$. In this case we can find integers $n\ge1$
and $k\ge0$ such that \[
\lambda_{n-1}\le t\le\lambda_{n}\le\lambda_{n+k}\le t'\le\lambda_{n+k+1}\,.\]
If $k\ge1$ then we have \begin{equation}
0<g(x,t')-g(x,t)=\left[g(x,t')-g(x,\lambda_{n+k})\right]+\sum_{m=1}^{k}\left[g(x,\lambda_{n+m})-g(x,\lambda_{n+m-1})\right]+\left[g(x,\lambda_{n})-g(x,t)\right]\,.\label{eq:bwt}\end{equation}
If $k=0$ then we have the same equation but with the middle sum $\sum_{m=1}^{k}\left[g(x,\lambda_{n+m})-g(x,\lambda_{n+m-1})\right]$
deleted. We can apply the preceding first step of this proof separately
to each term in square brackets on the right side of (\ref{eq:bwt})
to show that the whole right side is dominated by

\[
\left[t'-\lambda_{n+k}\right]+\sum_{m=1}^{k}\left[\lambda_{n+m}-\lambda_{n+m-1}\right]+\left[\lambda_{n}-t\right]\,,\]
where again the middle sum is deleted if $k=0$. Since this last expression
equals $t'-t$, our proof of (\ref{eq:xlip}), and therefore also
of (\ref{eq:slipc}), has now been completed.

We can finally give the definition of the operator $S_{N}$ which
acts on the one dimensional space $\left\{ c\chi_{I_{N}}:c\in\mathbb{C}\right\} $.
For each complex constant $c$, we have that $S_{N}\left(c\chi_{I_{N}}\right)$
is the restriction to the interval $[0,1]$ of the function \[
x\mapsto g\left(x,\gamma\left(\left|c\right|w\right)\right)\,.\]

\begin{rem}
We may care to remember that the functions $g$ and $\gamma$ used
here both depend crucially on the sequence $\left\{ h_{n}\right\} $
and the other sequences derived from it. Therefore they depend on
the the number $w=2^{-Np}$. Nevertheless, we shall establish some
very useful estimates and properties of $S_{N}$ which do not depend
on $N$. For example, in Lemma \ref{lem:lip}, we will benefit from
the fact that, unlike the above mentioned sequences, the sequence
$\left\{ m_{n}\right\} _{n\ge0}$ does not depend on $w$.
\end{rem}
Since $\gamma$ is an increasing function and $t\mapsto g(x,t)$ is
a non decreasing function of $t$ for each fixed $x$, we immediately
obtain that $S_{N}$ has the pointwise monotonicity property that\begin{equation}
S_{N}(c\chi_{I_{N}})\le S_{N}(c'\chi_{I_{N}})\mbox{ whenever }\left|c\right|\le\left|c'\right|\,.\label{eq:monp}\end{equation}

Now we shall obtain an $L^{1}$ Lipschitz estimate for $S_{N}$, which
again uses the monotonicity of $t\mapsto g(x,t)$, and also (\ref{eq:mbet}).
Let $c$ and $c'$ be any two complex numbers. We can suppose, without
loss of generality, that $\left|c\right|\le\left|c'\right|$. Then
we have $g\left(x,\gamma\left(\left|c'\right|w\right)\right)-g\left(x,\gamma\left(\left|c\right|w\right)\right)\ge0$
and so a series of steps, using various properties of $g$, $G$ and
$\gamma$, including (\ref{eq:fgg}), will give us that \begin{eqnarray}
\left\Vert S_{N}(c'\chi_{I_{N}})-S_{N}(c\chi_{I_{N}})\right\Vert _{L^{1}} & = & \int_{0}^{1}\left|g\left(x,\gamma\left(\left|c'\right|w\right)\right)-g\left(x,\gamma\left(\left|c\right|w\right)\right)\right|dx\nonumber \\
 & \le & \int_{0}^{w+h_{0}}\left|g\left(x,\gamma\left(\left|c'\right|w\right)\right)-g\left(x,\gamma\left(\left|c\right|w\right)\right)\right|dx\nonumber \\
 & = & \int_{0}^{w+h_{0}}g\left(x,\gamma\left(\left|c'\right|w\right)\right)-g\left(x,\gamma\left(\left|c\right|w\right)\right)dx\nonumber \\
 & = & \left|G\left(\gamma\left(\left|c'\right|w\right)\right)\right|-\left|G\left(\gamma\left(\left|c\right|w\right)\right)\right|\nonumber \\
 & = & \left|c'\right|w-\left|c\right|w\le\left|c'-c\right|w=\left\Vert c'\chi_{I_{N}}-c\chi_{I_{N}}\right\Vert _{L^{1}}\,.\label{eq:pmpm}\end{eqnarray}

The $L^{\infty}$ boundedness of $S_{N}$ is also straightforward.
For each complex number $c$ we have, using simple properties of $g$
and (\ref{eq:mma}), that \begin{equation}
\left\Vert S_{N}\left(c\chi_{I_{N}}\right)\right\Vert _{L^{\infty}}\le\sup_{x\ge0}g\left(x,\gamma\left(\left|c\right|w\right)\right)=g\left(0,\gamma\left(\left|c\right|w\right)\right)=\gamma\left(\left|c\right|w\right)\le\left|c\right|=\left\Vert c\chi_{I_{N}}\right\Vert _{L^{\infty}}\,.\label{eq:lefs}\end{equation}

We can also obtain an $L^{\infty}$ Lipschitz estimate for $S_{N}$.
Here again we consider any two complex numbers $c$ and $c'$ and
we will proceed, using (\ref{eq:slipc}) and then (\ref{eq:gamlip}).
We see that \begin{eqnarray}
\left\Vert S_{N}(c'\chi_{I_{N}})-S_{N}(c\chi_{I_{N}})\right\Vert _{L^{\infty}} & \le & \sup_{x\ge0}\left|g(x,\gamma\left(\left|c'\right|w\right))-g(x,\gamma\left(\left|c\right|w\right))\right|\nonumber \\
 & \le & \left|\gamma\left(\left|c'\right|w\right))-\gamma\left(\left|c\right|w\right))\right|\nonumber \\
 & \le & \frac{1}{w}\left|\left|c'\right|w-\left|c\right|w\right|=\left|\left|c'\right|-\left|c\right|\right|\nonumber \\
 & \le & \left|c'-c\right|=\left\Vert c'\chi_{I_{N}}-c\chi_{I_{N}}\right\Vert _{L^{\infty}}\,.\label{eq:epm}\end{eqnarray}
\smallskip{}

The following result will help us later to establish that the operator
$T_{5}$ maps bounded subsets of $L^{\infty}$ into compact subsets
of $L^{\infty}$.
\begin{lem}
\label{lem:lip}For each positive constant $C$ there exists another
positive constant $L=L(C,p)$ depending only on $C$ and $p$, such
that, for all complex numbers $\alpha$ with $\left|\alpha\right|\le C$,
the function $S_{N}\left(\alpha\chi_{I_{N}}\right)$ satisfies a Lipschitz
condition with Lipschitz constant not exceeding $L(C,p)$.
\end{lem}
\textit{Proof.} Given $C$, let $n=n_{C}$ be the smallest positive
integer which satisfies $2^{n}-1\ge C$. (More explicitly, we have
$n_{C}=\left\lceil \frac{\log(C+1)}{\log2}\right\rceil $~.) Then,
for each $\alpha$ which satisfies $\left|\alpha\right|\le C$, we
use the monotonicity of $\gamma$ and (\ref{eq:twsu}) to obtain that
\begin{equation}
\gamma(\left|\alpha\right|w)\le\gamma\left((2^{n_{C}}-1)w\right)=\lambda_{n_{C}}\,.\label{eq:jsat}\end{equation}
Let us fix $t=\gamma(\left|\alpha\right|w)$. Then, in view of (\ref{eq:jsat}),
the integer $\nu(t)$, which is defined as in (\ref{eq:defnu}), must
satisfy $\nu(t)\le n_{C}$. Therefore, since the sequence $\left\{ m_{n}\right\} _{n\ge0}$
is increasing, we have $m_{\nu(t)}\le m_{n_{C}}$. We combine this
with (\ref{eq:lipc}) to obtain that the function $x\mapsto g(x,t)$
satisfies a Lipschitz condition on $[0,\infty)$ with Lipschitz constant
not exceeding $m_{n_{C}}$. Since $S_{N}\left(\alpha\chi_{I_{N}}\right)$
is the restriction of this function to $[0,1]$, our proof is complete,
with the constant $L(C,p)$ given by \begin{equation}
L(C,p)=m_{n_{C}}=\sum_{k=1}^{n_{C}}2^{k(p+1)}=\frac{2^{\left(\left\lceil \frac{\log(C+1)}{\log2}\right\rceil +1\right)(p+1)}-2^{(p+1)}}{2^{p+1}-1}\,.\label{eq:dlcp}\end{equation}
$\qed$

\bigskip{}

Recalling that $w=2^{-Np}$ we observe that $h_{N-1}=\sqrt{w^{2}+\frac{w}{2^{pN}}}-w=\sqrt{w^{2}+w^{2}}-w=(\sqrt{2}-1)w$.
Therefore \begin{equation}
y_{N}=\frac{2^{N-1}}{1+\frac{h_{N-1}}{2w}}=\frac{2^{N-1}}{1+\frac{\sqrt{2}-1}{2}}=\frac{2^{N}}{1+\sqrt{2}}\,.\label{eq:yn}\end{equation}
We will need the preceding formula for our next step. This will be
to consider the particular function $(2^{N}-1)\chi_{I_{N}}$ which
of course satisfies \begin{equation}
\left\Vert (2^{N}-1)\chi_{I_{N}}\right\Vert _{L^{p}}\le\left\Vert 2^{N}\chi_{I_{N}}\right\Vert _{L^{p}}=1\,.\label{eq:nib}\end{equation}
In preparation for showing later that the operator $T$ does not map
bounded subsets of $L^{p}$ into compact subsets of $L^{p}$, we shall
estimate the norm $\left\Vert S_{N}\left((2^{N}-1)\chi_{I_{N}}\right)\right\Vert _{L^{p}}$
from below. With the help of (\ref{eq:twsu}), the definitions of
the function $g$ and the sequence $\left\{ \lambda_{n}\right\} _{n\ge0}$
and then finally (\ref{eq:yn}), we see that, for all points $x$
in the interval $[0,w]=[0,2^{-Np}]$, the function $g\left(x,\gamma\left(\left(2^{N}-1\right)w\right)\right)$
satisfies \begin{eqnarray*}
g\left(x,\gamma\left(\left(2^{N}-1\right)w\right)\right) & = & g\left(x,\lambda_{N}\right)=\lambda_{N}\ge y_{N}=\frac{2^{N}}{1+\sqrt{2}}\,.\end{eqnarray*}
This means that \begin{eqnarray*}
\left\Vert S_{N}\left((2^{N}-1)\chi_{I_{N}}\right)\right\Vert _{L^{p}}^{p} & \ge & \left\Vert S_{N}\left((2^{N}-1)\chi_{I_{N}}\right)\cdot\chi_{[0,w]}\right\Vert _{L^{p}}^{p}\\
 & = & \int_{0}^{w}g\left(x,\gamma\left(\left(2^{N}-1\right)w\right)\right)^{p}dx\\
 & \ge & \frac{2^{Np}w}{\left(1+\sqrt{2}\right)^{p}}=\frac{1}{\left(1+\sqrt{2}\right)^{p}}\end{eqnarray*}
and we have shown that \begin{equation}
\left\Vert S_{N}\left((2^{N}-1)\chi_{I_{N}}\right)\right\Vert _{L^{p}}\ge\frac{1}{1+\sqrt{2}}\ \mbox{for each }N\in\mathbb{N}\,.\label{eq:ntl}\end{equation}

\subsubsection{\label{sub:together}Putting all the pieces together}

Now that we have constructed and described the properties of the special
operators $S_{N}$ we can turn to showing that the operator $T_{5}$
obtained from those operators by the formula (\ref{eq:deft}) has
all the properties needed to make it the counterexample that we are
seeking. In this subsection we will often simply write $T$ instead
of $T_{5}$.

First we consider the action of $T$ on the zero function. Since $S_{N}(0)=0$
for each $N$ we deduce that $T(0)=0$.

Next we observe that, for any for any two functions $f$ and $g$
in $L^{1}$ which satisfy $\left|f(x)\right|\le\left|g(x)\right|$
for almost every $x$, we of course have $0\le\int_{I_{N}}\left|f(x)\right|dx\le\int_{I_{N}}\left|g(x)\right|dx$
for each $N$. So, with the help of (\ref{eq:monp}), we obtain the
pointwise estimate \begin{equation}
0\le T(f)\le T(g)\ \mbox{whenever }\left|f(x)\right|\le\left|g(x)\right|\,\mbox{for a.e. }x\in(0,1)\,.\label{eq:smonp}\end{equation}

This property will now help us show that $T$ satisfies Lipschitz
norm estimates for both $L^{1}$ and $L^{\infty}$, namely that \begin{equation}
\left\Vert Tf_{1}-Tf_{2}\right\Vert _{L^{q}}\le\left\Vert f_{1}-f_{2}\right\Vert _{L^{q}}\,\mbox{for all }f_{1},f_{2}\in L^{q}\,\mbox{and for }q=1,\infty\,.\label{eq:lipe}\end{equation}

For each such $f_{1}$ and $f_{2}$ and $q$ we obviously have $\left\Vert \left|f_{1}\right|-\left|f_{2}\right|\right\Vert _{L^{q}}\le\left\Vert f_{1}-f_{2}\right\Vert _{L^{q}}$.
Furthermore, $Tf=T\left(\left|f\right|\right)$ for each $f\in L^{1}$.
This means that it suffices to prove (\ref{eq:lipe}) in the special
case where $f_{1}$ and $f_{2}$ are both non negative functions.
For two such functions let us set $f_{-}=\min\left\{ f_{1},f_{2}\right\} $
and $f_{+}=\max\left\{ f_{1},f_{2}\right\} $. Then $\left|f_{1}(x)-f_{2}(x)\right|=f_{+}(x)-f_{-}(x)$
and also, by (\ref{eq:smonp}), we have the two pointwise estimates
$T(f_{-})\le T(f_{j})\le T(f_{+})$ for $j=1,2$ which imply that
$\left|T(f_{1})-T(f_{2})\right|\le T(f_{+})-T(f_{-})$. From all this
we see that it will suffice to prove (\ref{eq:lipe}) in the special
case where $0\le f_{1}\le f_{2}$.

Let $Q$ be the (linear) conditional expectation operator defined
by \begin{equation}
Qf=\sum_{N=1}^{\infty}\frac{1}{2^{-Np}}\int_{I_{N}}f(x)dx\cdot\chi_{I_{N}}\,\mbox{ for each }f\in L^{1}\,.\label{eq:defq}\end{equation}
Obviously $\left\Vert Qf_{1}-Qf_{2}\right\Vert _{L^{q}}\le\left\Vert f_{1}-f_{2}\right\Vert _{L^{q}}$
for all $f_{1},f_{2}\in L^{q}$ when $q=1$ and when $q=\infty$.
Furthermore $Tf=T\left(Qf\right)$ for all non negative $f\in L^{1}$.
This enables us to further reduce the proof of (\ref{eq:lipe}) to
a still more special case. Not only does it suffice to consider $f_{1}$
and $f_{2}$ satisfying $0\le f_{1}\le f_{2}$. We can also suppose
that $f$$_{1}$ and $f_{2}$ are both functions of the form $\sum_{N=1}^{\infty}\alpha_{N}\chi_{I_{N}}$.

Let $\epsilon$ be an arbitrary positive number. For each $N\in\mathbb{N}$
let $H_{N}$ be the measurable set \[
H_{N}=\left\{ x\in[0,1]:S_{N}(f_{2}\chi_{I_{N}})(x)\ge T_{5}(f_{2})(x)-\epsilon\right\} \,.\]
Since $T(f_{2})(x)=\sup_{N\in\mathbb{N}}S(f_{2}\chi_{I_{N}})(x)$,
we have $\bigcup_{N\in\mathbb{N}}H_{N}=[0,1]$. Now we use the sequence
of sets $\left\{ H_{N}\right\} _{N\in\mathbb{N}}$ to obtain another
sequence $\left\{ \Omega_{N}\right\} _{N\in\mathbb{N}}$ of pairwise
disjoint measurable sets such that $\Omega_{N}\subset H_{N}$ for
each $N$ and $\bigcup_{N\in\mathbb{N}}\Omega_{N}=[0,1]$. We can
do this in the usual and obvious way, by setting $\Omega_{1}=H_{1}$
and then proceeding recursively by taking $\Omega_{N}=H_{N}\backslash\left(\bigcup_{k=1}^{N-1}\Omega_{k}\right)$
for each $N\ge2$. (Of course some of the sets $\Omega_{N}$ may be
empty.)

For each $N\in\mathbb{N}$ and for each $x\in\Omega_{N}$ we have
\[
0\le T(f_{2})(x)T(f_{1})(x)\le T(f_{2})(x)-S_{N}(f_{1}\chi_{I_{N}})(x)\le\epsilon+S_{N}(f_{2}\chi_{I_{N}})(x)-S_{N}(f_{1}\chi_{I_{N}})(x)\,.\]
This means that, for $q=1,\infty$, we have \begin{eqnarray*}
\left\Vert \left(T(f_{2})-T(f_{1})\right)\cdot\chi_{\Omega_{N}}\right\Vert _{L^{q}} & \le & \left\Vert \epsilon\chi_{\Omega_{N}}+S_{N}(f_{2}\chi_{I_{N}})-S_{N}(f_{1}\chi_{I_{N}})\right\Vert _{L^{q}}\\
 & \le & \epsilon\left\Vert \chi_{\Omega_{N}}\right\Vert _{L^{q}}+\left\Vert S_{N}(f_{2}\chi_{I_{N}})-S_{N}(f_{1}\chi_{I_{N}})\right\Vert _{L^{q}}\,.\end{eqnarray*}
Now we can apply the Lipschitz norm estimates (\ref{eq:pmpm}) if
$q=1$ or (\ref{eq:epm}) if $q=\infty$, to obtain that \begin{equation}
\left\Vert \left(T(f_{2})-T(f_{1})\right)\chi_{\Omega_{N}}\right\Vert _{L^{q}}\le\epsilon\left\Vert \chi_{\Omega_{N}}\right\Vert _{L^{q}}+\left\Vert (f_{2}\chi_{I_{N}})-(f_{1}\chi_{I_{N}})\right\Vert _{L^{q}}\,.\label{eq:ctt}\end{equation}
In the case where $q=1$ we sum both sides of the preceding inequality
over all $N$ and obtain that \begin{eqnarray*}
\left\Vert T(f_{2})-T(f_{1})\right\Vert _{L^{1}} & = & \sum_{N=1}^{\infty}\left\Vert \left(T(f_{2})-T(f_{1})\right)\chi_{\Omega_{N}}\right\Vert _{L^{1}}\\
 & \le & \epsilon\sum_{N=1}^{\infty}\left\Vert \chi_{\Omega_{N}}\right\Vert _{L^{1}}+\sum_{N=1}^{\infty}\left\Vert (f_{2}\chi_{I_{N}})-(f_{1}\chi_{I_{N}})\right\Vert _{L^{1}}\\
 & = & \epsilon+\left\Vert \left(f_{2}-f_{1}\right)\chi_{\bigcup_{N\in\mathbb{N}}I_{N}}\right\Vert _{L^{1}}=\epsilon+\left\Vert f_{2}-f_{1}\right\Vert _{L^{1}}\,.\end{eqnarray*}
In the case where $q=\infty$ we take the supremum over all $N$ of
both sides of (\ref{eq:ctt}) to obtain that\begin{eqnarray*}
\left\Vert T(f_{2})-T(f_{1})\right\Vert _{L^{\infty}} & = & \sup_{N\in\infty}\left\Vert \left(T(f_{2})-T(f_{1})\right)\chi_{\Omega_{N}}\right\Vert _{L^{\infty}}\\
 & \le & \epsilon\sup_{N\in\infty}\left\Vert \chi_{\Omega_{N}}\right\Vert _{L^{\infty}}+\sup_{N\in\infty}\left\Vert (f_{2}\chi_{I_{N}})-(f_{1}\chi_{I_{N}})\right\Vert _{L^{\infty}}\\
 & = & \epsilon+\left\Vert \left(f_{2}-f_{1}\right)\chi_{\bigcup_{N\in\mathbb{N}}I_{N}}\right\Vert _{L^{\infty}}=\epsilon+\left\Vert f_{2}-f_{1}\right\Vert _{L^{\infty}}\,.\end{eqnarray*}
Since we may take $\epsilon$ to be arbitrarily small, the preceding
calculations establish (\ref{eq:lipe}) in the special case specified
above, which, as already explained, also suffices to prove (\ref{eq:lipe})
in full generality. Since $T(0)=0$ we also know from (\ref{eq:lipe})
that \begin{equation}
\left\Vert T(f)\right\Vert _{L^{q}}\le\left\Vert f\right\Vert _{L^{q}}\,\mbox{for all }f\in L^{q}\,\mbox{and for }q=1,\infty\,.\label{eq:bb}\end{equation}

We need to show that our operator $T$ maps bounded subsets of $L^{\infty}$
into relatively compact subsets of $L^{\infty}$. One ingredient for
doing that will be the following simple result. It is surely a special
case of well known and more general results. But it seems just as
easy to prove it as to give a reference.
\begin{lem}
\label{lem:loli}Let $L$ and $C$ be positive constants and let $\left\{ f_{N}\right\} _{N\in\mathbb{N}}$
be a sequence of functions $f_{N}:[0,1]\to\mathbb{\mathbb{R}}$ which
all satisfy $\left|f_{N}(x)\right|\le C$ and $\left|f_{N}(x)-f_{N}(x')\right|\le L\left|x-x'\right|$
for all $x,x'\in[0,1]$. Then the function $g:[0,1]\to\mathbb{R}$
defined by $g(x)=\sup_{N\in\mathbb{N}}f_{N}(x)$ also satisfies $\left|g(x)\right|\le C$
and $\left|g(x)-g(x')\right|\le L\left|x-x'\right|$ for all $x,x'\in[0,1]$.
\end{lem}
\textit{Proof.} As our first step, consider two arbitrary functions
$u_{1}:[0,1]\to\mathbb{R}$ and $u_{2}:[0,1]\to\mathbb{R}$ which
satisfy $\left|u_{j}(x)-u_{j}(x')\right|\le L\left|x-x'\right|$ for
$j=1,2$ and all $x,x'\in[0,1]$. Let $w=\max\left\{ u_{1},u_{2}\right\} $.
We shall show that \begin{equation}
\left|w(x)-w(x')\right|\le L\left|x-x'\right|\label{eq:wlip}\end{equation}
 for each $x,x'\in[0,1]$. We may suppose, without loss of generality,
that $x<x'$. If the continuous function $t\mapsto u_{1}(t)-u_{2}(t)$
has the same sign at both endpoints of the interval $[x,x']$ or vanishes
at one of these endpoints, then $w(x)-w(x')$ equals either $u_{1}(x)-u_{1}(x')$
or $u_{2}(x)-u_{2}(x')$ and in either of these cases we obtain (\ref{eq:wlip}).
Otherwise there must be some point $x"\in(x,x')$ for which $u_{1}(x")-u_{2}(x")=0$
and so we can apply the preceding argument on each of the intervals
$[x,x"]$ and $[x",x']$ to show that $\left|w(x)-w(x")\right|\le L\left|x-x"\right|=L\left(x"-x\right)$
and $\left|w(x')-w(x")\right|\le L\left|x'-x"\right|=L(x'-x")$, which
together imply (\ref{eq:wlip}).

For our second and final step we observe that, by simply reiterating
the previous step, we can obtain, for each $N\in\mathbb{N}$, that
the function $g_{N}=\max\left\{ f_{1},f_{2},....,f_{N}\right\} $
satisfies $\left|g_{N}(x)-g_{N}(x')\right|\le L\left|x-x'\right|$
for all $x,x'\in(0,1)$, and obviously it also satisfies $\left|g_{N}(x)\right|\le C$
for all $x\in(0,1)$. Since $g(x)=\lim_{N\to\infty}g_{N}(x)<\infty$,
we can pass to the limit in the two preceding inequalities to obtain
the two required properties of $g$. $\qed$\medskip{}

Now suppose that $A$ is some bounded subset of $L^{\infty}$. Let
$C=\sup_{f\in A}\left\Vert f\right\Vert _{L^{\infty}}$. Let $B=\left\{ Q\left(\left|f\right|\right):f\in A\right\} $
where $Q$ is the conditional expectation operator defined above.
Then obviously $\sup_{f\in B}\left\Vert f\right\Vert _{L^{\infty}}\le C$
and $T(A)=T(B)$. For each $f\in B$ we can apply Lemma \ref{lem:lip}
to obtain that, for each $N\in\mathbb{N}$, the function $S_{N}\left(f\chi_{I_{N}}\right)$
satisfies a Lipschitz condition with Lipschitz constant not exceeding
the number $L(C,p)$ defined in (\ref{eq:dlcp}). We also have $\left\Vert S_{N}\left(f\chi_{I_{N}}\right)\right\Vert _{L^{\infty}}\le C$,
in view of (\ref{eq:lefs}). These two facts enable us to apply Lemma
\ref{lem:loli} to obtain that $T(f)$ is also bounded by $C$ and
satisfies a Lipschitz condition with constant not exceeding $L(C,p)$.
Thus we have shown that $T(B)$ is a bounded and equicontinuous subset
of the Banach space $C[0,1]$ of continuous real valued functions
on $[0,1]$ equipped with the supremum norm. Therefore, by the Arzelà--Ascoli
theorem, $T(B)$ is a relatively compact subset of $C[0,1]$, and
therefore also of $L^{\infty}[0,1]$.

We have now reached the very last part of our discussion of the operator
$T=T_{5}$. It remains only to show that it does not map every bounded
subset of $L^{p}$ into a relatively compact subset of $L^{p}$. For
this we consider the particular set $A$ consisting of all the functions
$\psi_{N}=\left(2^{N}-1\right)\chi_{I_{N}}$ for all $N\in\mathbb{N}$.
We have already observed in (\ref{eq:nib}) that this is a bounded
subset of $L^{p}$. Since $T(\psi_{N})=S_{N}(\psi_{N})$, we also
know from (\ref{eq:ntl}) that $\left\Vert T(\psi_{N})\right\Vert _{L^{p}}\ge\frac{1}{1+\sqrt{2}}$
for each $N$. If $T\left(A\right)$ is relatively compact in $L^{p}$
then some subsequence $\left\{ T\left(\psi_{N_{k}}\right)\right\} _{k\in\mathbb{N}}$
must converge to some function $\phi$ in $L^{p}$ norm. In view of
Hölder's inequality, the same subsequence must converge to $\phi$
also in $L^{1}$ norm. Since \[
\left\Vert T\left(\psi_{N}\right)\right\Vert _{L^{1}}\le\left\Vert \psi_{N}\right\Vert _{L^{1}}=(2^{N}-1)2^{-Np}\le2^{-N(p-1)},\]
the function $\phi$ must be the zero function. But the above mentioned
strictly positive bound from below for $\left\Vert T\left(\psi_{N}\right)\right\Vert _{L^{p}}$
means that $\left\{ T\left(\psi_{N_{k}}\right)\right\} _{k\in\mathbb{N}}$
cannot converge to $0$ in $L^{p}$ norm. This proves that the set
$T(A)$, i.e., the set $T_{5}(A)$, cannot be relatively compact in
$L^{p}$.

\subsubsection{A modification of this example showing that two sided compactness
conditions are also insufficient.}

The reader who has kept us company till now, may be interested to
know that it is possible to compose the operator $T_{5}$ which we
have just constructed, with another Lipschitz operator, so that the
new composed operator satisfies an additional compactness condition
at the other {}``endpoint'', i.e., property {[}iii$_{1}${]}, and
it still has all the other properties of $T_{5}$, namely {[}i{]},
{[}ii{]}, {[}iii$_{0}${]} and non compactness on $L^{p}$ for the
value of $p$, that we chose in advance.

We proceed somewhat analogously to the arguments used in Section \ref{sec:twosided}
to combine the examples of Sections \ref{sec:onesided} and \ref{sec:smallonesided}.

Having chosen our $p\in(1,\infty)$, we begin by constructing exactly
the same operator $T_{5}$ for that value of $p$ as was constructed
in the preceding subsections. Our new operator $T_{6}$, which will
have all the properties listed just above, will be the composition
$T_{6}=T_{5}\circ V$ of $T_{5}$ with another operator $V$ which
will be rather similar to the the operator in Section \ref{sec:onesided}.
But this time we let $\left\{ I_{N}\right\} _{N\in\mathbb{N}}$ denote
exactly that sequence of pairwise disjoint open subintervals of $(0,1)$
with $\left|I_{N}\right|=2^{-Np}$ which was introduced at the beginning
of the construction in Subsection \ref{sub:elaborate}. Let $Q$ be
the linear operator of conditional expectation with respect to this
sequence, as defined in (\ref{eq:defq}). Let $v:(0,1)\to[0,\infty)$
be the function $v=\sum_{N=1}^{\infty}(2^{N}-1)\chi_{I_{N}}$. Now
we can define the nonlinear operator $V$ by \[
V(f)=\min\left\{ \left|Qf\right|,v\right\} \,\mbox{for all }f\in L^{1}\,.\]

This time we let $H$ be the set of all functions $f:(0,1)\to[0,\infty)$
of the form $f=\sum_{n=1}^{\infty}\alpha_{n}\chi_{I_{n}}$ where each
of the constants $\alpha_{n}$ satisfies $0\le\alpha_{n}\le(2^{n}-1)$.
This time we can use the convergence of the series $\sum_{n=1}^{\infty}(2^{n}-1)\left|I_{n}\right|$
to show that $H$ is a compact subset of $L^{1}$. We know from our
previous discussion that $T_{5}$ satisfies $\left\Vert T_{5}(f)-T_{5}(g)\right\Vert _{L^{1}}\le\left\Vert f-g\right\Vert _{L^{1}}$.
So we deduce that the set $T_{5}(H)$, as the continuous image of
a compact set, is also a compact subset of $L^{1}$. Since $V(L^{1})=H$
we obtain that $T_{6}=T_{5}\circ V$ maps $L^{1}$ and therefore also
every bounded subset of $L^{1}$ into the compact subset $T_{5}(H)$.

Let $J$ be an arbitrary bounded subset of $L^{\infty}$. Then of
course $V(J)$ is also a bounded subset of $L^{\infty}$ and so, again
using a property of $T_{5}$ established above, we have that $T_{5}\left(V(J)\right)=T_{6}(J)$
is a relatively compact subset of $L^{\infty}$.

We now know that $T_{6}$ satisfies properties {[}iii$_{0}${]} and
{[}iii$_{1}${]}. Property {[}i{]} is obvious and property {[}ii{]}
follows trivially from the fact that $T_{6}$ and $V$ both satisfy
{[}ii{]}.

Finally we show that the map $T_{6}:L^{p}\to L^{p}$ is not compact.
As in our previous treatment of $T_{5}$, we again consider the set
$A$ consisting of all the functions $\psi_{N}=\left(2^{N}-1\right)\chi_{I_{N}}$
for all $N\in\mathbb{N}$. We already know that this is a bounded
subset of $L^{p}$ and that its image $T_{5}(A)$ is not a relatively
compact subset of $L^{p}$. It remains to make the trivial observation
that $V(\psi_{N})=\psi_{N}$ for each $N$ and therefore $T_{6}(A)=T_{5}(A)$.~

\smallskip{}

\end{document}